\def\version{08/02/2005 Version 1}
\documentclass[11pt,a4paper]{amsart}
\usepackage{fullpage}
\usepackage{amssymb,amscd}
\usepackage[all]{xy}\CompileMatrices

\usepackage{mathrsfs}

\newtheorem{thm}{Theorem}[section]
\newtheorem{lem}[thm]{Lemma}
\newtheorem{prop}[thm]{Proposition}
\newtheorem{cor}[thm]{Corollary}

\theoremstyle{definition}
\newtheorem{rem}[thm]{Remark}
\newtheorem{defn}[thm]{Definition}
\newtheorem{ex}[thm]{Example}

\newtheorem*{cond*}{Condition}
\newenvironment{Cond}[1]{\par\medskip\noindent{\bf Condition #1.}
}{\par \medskip}

\newtheorem*{Quest*}{Question}
\newtheorem{assump}[thm]{Assumption}

\numberwithin{equation}{section}
\numberwithin{figure}{section}

\def\ie{\emph{i.e.}}
\def\ds{\displaystyle}
\def\:{\colon}
\def\.{\cdot}

\def\<{\left\langle}
\def\>{\right\rangle}
\def\({\left(}
\def\){\right)}
\def\ph#1{\phantom{#1}}
\def\epsilon{\varepsilon}
\def\phi{\varphi}
\def\leq{\leqslant}
\def\geq{\geqslant}
\def\lla{\longleftarrow}

\def\lra{\longrightarrow}
\def\Lra{\Longrightarrow}
\def\ra{\rightarrow}

\def\hat#1{\widehat{#1}}

\def\iso{\cong}

\DeclareMathOperator{\im}{im}

\def\F{\mathbb{F}}

\def\k{\Bbbk}

\def\Z{\mathbb{Z}}
\def\Smash_#1{\wedge_{#1}}

\def\Times_#1{\ds\mathop{\times}_{#1}}
\def\oTimes_#1{\otimes_{#1}}
\def\oPlus_#1{\bigoplus_{#1}}
\def\ideal{\triangleleft}

\DeclareMathOperator{\Hom}{Hom}

\DeclareMathOperator{\Pic}{Pic}

\DeclareMathOperator{\Tor}{Tor}
\def\id{\mathrm{id}}
\DeclareMathOperator*{\colim}{colim}

\DeclareMathOperator*{\holim}{holim}

\DeclareMathOperator{\cofibre}{cofibre}

\title[Invertible modules for some commutative $\mathbb{S}$-algebras]
{Invertible modules for commutative $\mathbb{S}$-algebras with
residue fields}
\author{Andrew Baker \and Birgit Richter}
\address{Mathematics Department, University of Glasgow,
Glasgow G12 8QW, Scotland.}
\email{a.baker@maths.gla.ac.uk}
\urladdr{http://www.maths.gla.ac.uk/$\sim$ajb}
\address{Mathematisches Institut der Universit\"at Bonn,
53115 Bonn, Germany.}
\email{richter@math.uni-bonn.de}
\urladdr{http://www.math.uni-bonn.de/people/richter}
\date{\version}
\thanks{We would like to thank Halvard Fausk, John Greenlees, Peter
Kropholler, and Stefan Schwede. We are grateful that John Rognes
spotted an embarrassing blunder in an earlier version of this paper.
The second author thanks the Mathematics Department of the University
of Glasgow and the University of Oslo for their hospitality; she was
supported by the \emph{Strategisk Universitetsprogram i Ren Matematikk}
(SUPREMA) of the Norwegian Research Council.}
\begin{document}
\begin{abstract}
The aim of this note is to understand under which conditions invertible
modules over a commutative $\mathbb{S}$-algebra in the sense of Elmendorf,
Kriz, Mandell \& May give rise to elements in the algebraic Picard group
of invertible graded modules over the coefficient ring by taking homotopy
groups. If a  connective commutative
$\mathbb{S}$-algebra $R$ has coherent localizations $(R_*)_{\mathfrak{m}}$
for every maximal ideal $\mathfrak{m} \ideal R_*$, then for every invertible
$R$-module~$U$, $U_*=\pi_*U$ is an invertible graded $R_*$-module. In
some non-connective cases we can carry the result over under the additional
assumption that the commutative $\mathbb{S}$-algebra has `residue fields' for 
all maximal ideals $\mathfrak{m}\ideal R_*$ if
the global dimension of $R_*$ is small or if $R$ is $2$-periodic with
underlying Noetherian complete local regular ring $R_0$.
\end{abstract}
\maketitle

\section*{Introduction}

For an arbitrary symmetric monoidal category $\mathscr C$, one can ask
which objects are invertible. The \emph{Picard group}, $\Pic(\mathscr C)$,
is then the collection  of isomorphism classes of such invertible objects
in $\mathscr C$. This does not have to be a set in general, but if it is
one, then $\Pic(\mathscr C)$ is an abelian group in a natural way.

The notion of Picard groups originates from algebraic geometry. The classical
example is that of the Picard groups of the  category of $A$-modules 
for a commutative ring~$A$. In recent years, topologists have introduced 
symmetric monoidal
categories of  spectra, the categories of modules over commutative
$\mathbb{S}$-algebras, whose homotopy theories are also symmetric monoidal
and provide natural generalizations of categories of modules over commutative
rings.

The paper by Strickland~\cite{NS:Picard}, following a talk of Hopkins,
introduced Picard groups in that framework. Examples have been considered
and in some  cases calculated in~\cite{NS:Picard,HMS,JPM:PicardGps,FLM}.

In general, it is \emph{not} clear if invertible modules over a commutative
$\mathbb{S}$-algebra in the sense of~\cite{EKMM} give rise to invertible
modules over the coefficient ring when one applies homotopy groups.
In~\ref{ex:KU} we give an explicit example for an invertible module over
a commutative $\mathbb{S}$-algebra $R$ whose coefficient module is not an
invertible graded $R_*$-module. We investigate some restricted classes of
commutative $\mathbb{S}$-algebras for which we can prove such a transfer
result. Similar questions about various other Picard groups were considered
in~\cite{Fausk,FLM,JPM:PicardGps}, but we are unaware of published topological
results of the form we describe. Our motivation comes in part from applications
to work on Galois extensions in~\cite{AB&BR:Galois}. However, our results apply
to classes of examples such as complex cobordism $MU$ and complex
$K$-theory.

We discuss mainly three classes of examples of commutative
$\mathbb{S}$-algebras for which we can prove a transfer result from
the Picard group of module spectra over these commutative
$\mathbb{S}$-algebras to the Picard group of graded modules over its
coefficients. Essentially these are as follows, 
\begin{itemize}
\item
connective commutative $\mathbb{S}$-algebras with coherent coefficient
rings or coefficient rings satisfying Eilenberg's condition,
\item
commutative $\mathbb{S}$-algebras with coherent coefficients, 
multiplicative residue fields and small global dimension, and
\item
commutative $\mathbb{S}$-algebras whose coefficients $R_*$ are of the
form $R_0[u,u^{-1}]$ with $R_0$ a Noetherian complete local regular ring.
\end{itemize}
For the precise statements see Theorems~\ref{thm:Invtble-Main}, 
\ref{thm:Eilenberg-R}, \ref{thm:nonconnMain}, and \ref{thm:RegNoethLocal}.


In the non-connective case we add the restriction that the commutative
$\mathbb{S}$-algebra $R$ has `residue fields' for all maximal ideals in its
coefficients. Although we do not have a complete understanding of such
$\mathbb{S}$-algebras, it appears that most standard examples satisfy this
requirement; however, in Example~\ref{non-ex} we draw attention to a spectrum
which has no such residue field. In Section~\ref{sec:Examples} we summarize our
results by providing a list of examples of commutative
$\mathbb{S}$-algebras for which the topological and algebraic Picard
groups agree. We close with an explicit counterexample, proving that
$\Pic(R)$ might differ from $\Pic(R_*)$ in general.


\section{The Picard group of a commutative $\mathbb{S}$-algebra}\label{sec:Pic}

In this section we relate the Picard group
of the commutative $\mathbb{S}$-algebra $R$ to the Picard group of
 the coefficient ring $R_*$. Of course this is a special
case of the more general notion for a symmetric monoidal
category~\cite{Fausk,FLM,JPM:PicardGps}. 

We follow~\cite{NS:Picard,HMS} in defining $\Pic(R)$ to be the collection
of equivalence classes of invertible objects in the derived category of 
$R$-modules, $\mathscr{D}_{R}$. By Proposition~\ref{prop:StrongDual-Cell2}, 
$\Pic(R)$ is a set 
in all cases that we will consider. Whenever choosing a representative for an 
equivalence class $[U] \in \Pic(R)$ we will pick a cofibrant $R$-module.  
Defining the product of equivalence classes $[U]$ and $[V]$ by
\[
[U][V]=[U\wedge_R V],
\]
$\Pic(R)$ becomes an abelian group. We also have an algebraic Picard
group $\Pic(R_*)$ of invertible graded $R_*$-modules and our goal is
to discuss the relationship of this with $\Pic(R)$.

Let $M_*$ be an invertible graded $R_*$-module, say with inverse $N_*$,
so
\[
M_*\otimes_{R_*}N_*\iso R_*.
\]
Then $M_*$ and $N_*$ are finitely
generated projective of constant rank~$1$. Choosing a finitely generated
free cover $F_*\lra M_*$, we see that $M_*$ is a summand of $F_*$. Of
course we have $F_*=\pi_*F$, where $F$ is a finite wedge of suspensions
of copies of $R$-spheres $S_R$. Furthermore, the associated 
splitting can be realized by a homotopy idempotent self-map $\epsilon\:F\lra 
F$. Now $M_*$ is realized as the homotopy image of $\epsilon$. (This argument 
was pointed out to us by the referee.)

Given the above realization of $M_*$ as an $R$-module $M$, we obtain a map
\begin{equation}\label{eqn:PicAlg->Top}
\Phi\:\Pic(R_*)\lra\Pic(R);\quad \Phi([M_*])=[M]
\end{equation}
which is a group homomorphism. If $[M_*]\in\ker\Phi=[R_*]$, then $M\simeq R$
and so $[M_*]=[R_*]=0$. Thus $\Phi$ is a monomorphism. The main aim of this
paper is to establish conditions under which $\Phi$ is an
isomorphism.

\section{Recollections on coherent rings and modules}\label{sec:Coherent}

The notion of a \emph{coherent commutative ring} has proved important in
topology, especially in connection with $MU$. We begin by reviewing the
basic notions. As a convenient reference we cite Cohen~\cite{JCohen:Coherent},
but the algebraic theory can be found in many places such
as~\cite{Bourbaki:CommAlg,Serre:CoherentSheaves}.

Let $A$ be a (possibly graded) commutative ring. Then an $A$-module $M$
is \emph{finitely presented} if there is a short exact sequence
\[
0\ra K\lra F\lra M\ra0
\]
in which $F$ and $K$ are finitely generated and $F$ is free. Such
a short exact sequence is called a \emph{finite presentation of $M$}.
\begin{lem}\label{lem:Coherent-ShanuelLemma}
Let $M$ be finitely presented and suppose that
\[
0\ra L\lra P\lra M\ra0
\]
is a short exact sequence in which $P$ is finitely generated and
projective. Then $L$ is finitely generated.
\end{lem}
\begin{proof}
Let
\[
0\ra K\lra F\lra M\ra0
\]
be a finite presentation of $M$. By Schanuel's Lemma, there is an
isomorphism
\[
P\oplus K\iso F\oplus L.
\]
Now since the left-hand side is finitely generated, $L$ is also
finitely generated.
\end{proof}

The $A$-module $M$ is \emph{coherent} if it and all its finitely
generated submodules are finitely presented. $A$ is \emph{coherent}
if it is coherent as an $A$-module.

Let $A$ be a commutative local ring and $M$ an $A$-module. Recall
that a resolution $F_* \lra M$ is \emph{minimal} if  for each~$n$,
the differential $d_n\:F_{n}\lra F_{{n-1}}$ has
$\ker d_n\subseteq\mathfrak{m}F_{n}$ and so $\im d_n\subseteq\mathfrak{m}F_{{n-1}}$.
\begin{prop}\label{prop:Coherent-MinResn}
Let $A$ be a commutative local ring with maximal ideal $\mathfrak m$.
If $M$ is a coherent $A$-module, then $M$ admits a minimal resolution
$F_*\lra M\ra0$ which is by finitely generated free modules $F_n$.
\end{prop}
\begin{proof}
We begin by choosing a finitely generated free module $F_0$ with the
property that its reduction modulo $\mathfrak{m}$ is isomorphic to
the reduction of the module $M$ modulo $\mathfrak{m}$; furthermore,
this isomorphism factors through $M$, giving the following diagram.
\[
\xymatrix{
{F_0} \ar[d]_{p} \ar[dr]\ar[r] & {F_0 \otimes A/\mathfrak{m}} \ar[d]^{\cong} \\
{M} \ar[r] & {M \otimes A/\mathfrak{m}}
}
\]
Note that by Nakayama's Lemma, the map $p$ is an epimorphism. As $p$
is a map between coherent modules, its kernel is finitely generated
and coherent. It is obvious that the  kernel of $p$ is contained in
$\mathfrak{m}F_0$. Following this pattern of  argument we can inductively
produce a minimal resolution as required.
\end{proof}

In our work we will make use of the following result
of~\cite[proposition~1.5]{JCohen:Coherent}.
\begin{prop}\label{prop:Coherent-colimits}
Let $A_\alpha$ be a filtered direct system of coherent commutative
rings such that $A=\colim_\alpha A_\alpha$ is flat over each $A_\alpha$.
Then $A$ is coherent.
\end{prop}
\begin{cor}\label{cor:Coherent-localn}
Let $A$ be a coherent commutative ring and $\Sigma$ a multiplicative
subset. Then the localization $A[\Sigma^{-1}]$ is a coherent ring.
In particular, for every prime ideal $\mathfrak{p}\ideal A$, the
localization $A_\mathfrak{p}$ is coherent.
\end{cor}
\begin{proof}
This follows from Proposition~\ref{prop:Coherent-colimits} since
such localizations are filtered direct colimits and are exact.
\end{proof}
\begin{ex}\label{ex:Coherent}
\begin{enumerate}
\item[]
\item
Any commutative Noetherian ring is coherent.
\item
Any countably generated polynomial ring over a Noetherian ring
is coherent since it is a colimit as in
Proposition~\ref{prop:Coherent-colimits}.
\item
In particular, $MU_*$ and $BP_*$ are coherent and so are all their
localizations and quotients with respect to finitely generated ideals.
\end{enumerate}
\end{ex}

\section{Local reductions of $R$-modules}\label{sec:New}


Consider the following condition on a maximal ideal $\mathfrak{m}\ideal R_*$.
\begin{Cond}{(A)}
There is an $R$-module $W$ for which the $R_*$-module $W_*=\pi_*W$ is
isomorphic to the residue field $R_*/\mathfrak{m}$. If such a $W$  exists, 
we will refer to it as a residue field. We will often choose one and denote it 
by $\kappa(\mathfrak{m})$. 
\end{Cond}


Notice that $\kappa(\mathfrak{m})$ is clearly $R_{\mathfrak{m}}$-local,
where $R_{\mathfrak{m}}$ denotes the commutative $R$-algebra associated
with the algebraic localization
\[
\pi_*(\ph{X})_{\mathfrak{m}}=(R_*)_{\mathfrak{m}}\oTimes_{R_*}\pi_*(\ph{X})
\]
of the  homotopy functor on $\mathscr{M}_R$, for details see~\cite{EKMM}.
So the existence of such an $R$-module is equivalent to the existence of
a corresponding  $R_{\mathfrak{m}}$-module.

We will be interested in $\mathbb{S}$-algebras $R$ for which Condition~(A)
is satisfied by \emph{all} maximal ideals $\mathfrak{m}\ideal R_*$. Here
are some examples.
\begin{ex}\label{ex:ResFds-1}
If $R$ is connective then each maximal ideal $\mathfrak{m}\ideal R_*$
has for its residue field a quotient field
$\k(\mathfrak{m})=R_0/\mathfrak{m}_0$ of $R_0$. There is a corresponding
Eilenberg-Mac~Lane $R$-algebra $H\k(\mathfrak{m})$ which we may take for
$\kappa(\mathfrak{m})$.
\end{ex}
\begin{ex}\label{ex:ResFds-2}
Let $p>0$ be a prime and $R=MU_{I_{p,n}}$, where
$I_{p,n} = (p,v_1,\ldots,v_{n-1})$ is the $n$-th invariant ideal for~$p$.
Then the Morava $K$-theory $K(n)$ with $1\leq n<\infty$ is such a spectrum
$\kappa(I_{p,n})$.
\end{ex}
\begin{ex}\label{ex:Invtble-KO[1/2]}
$KO[1/2]$ satisfies Condition~(A) for every maximal ideal in
\[
KO[1/2]_*=\Z[1/2,y,y^{-1}].
\]
The only maximal ideal containing an odd prime $p$ is $(p)$ and we may
take $\kappa(p)=KO\wedge M(p)$ where $M(p)$ is the usual mod~$p$ Moore
spectrum.
\end{ex}
\begin{ex}\label{non-ex}
To see a non-example, we refer to \cite[example 7.6]{BKS}: the topological
significance of this is that the Tate spectrum $\hat{\mathbb{H}}(BC_3,\F_3)$
is a commutative $\mathbb{S}$-algebra, but the maximal ideal generated by
the exterior generator in degree one does not give rise to a residue field.
We learned of this example from John Greenlees.
\end{ex}

Example~\ref{ex:Invtble-KO[1/2]} leads us to introduce another condition
on a commutative $\mathbb{S}$-algebra $R$ and a maximal ideal
$\mathfrak{m}\ideal R_*$ that turns out to be useful in our work.
\begin{Cond}{(B)}
There is an $R$-algebra $R'$ where $R'_*$ is a local ring with maximal
ideal $\mathfrak{m}'\ideal R'_*$ which satisfies Condition~(A) and whose
unit induces a local homomorphism $(R_*)_{\mathfrak{m}}\lra R'_*$ which
makes $R'_*$ a flat $(R_*)_{\mathfrak{m}}$-module.
\end{Cond}

If $\mathfrak{m}$ satisfies Condition~(A), then the localization map
$R\lra R_{\mathfrak{m}}$ satisfies Condition~(B).
\begin{rem}\label{rem:V(m)-Extension}
Suppose that $B$ is a commutative ring and $A\subseteq B$ is a subring
so that~$B$ is a finite $A$-algebra. It is standard that if $A$ is
Noetherian then so is~$B$. Conversely, if $B$ is Noetherian then $A$
is Noetherian by the Eakin-Nagata theorem~\cite[theorem~3.7]{Matsumura}.
In these cases $B$ is automatically flat over $A$.
\end{rem}
\begin{ex}\label{ex:Invtble-KO}
Consider $KO$ and the maximal ideals
\[
\mathfrak{m}\ideal KO_*=\Z[\eta,y,w,w^{-1}]/(2\eta,\eta^3,\eta y,y^2-4w).
\]
If $p$ is an odd prime, then the only maximal ideal containing~$p$ is~$(p)$
and as in Example~\ref{ex:Invtble-KO[1/2]} we may take $\kappa(p)=KO\wedge M(p)$.
The only maximal ideal containing~$2$ is $(2,\eta,y)$ and we may take the
obvious map $KO\lra KU_{(2)}$ and set $\kappa(2,\eta,y)=KU\wedge M(2)$.

Thus every maximal ideal containing an odd prime $\mathfrak{m}\ideal KO_*$
satisfies Condition~(B).
\end{ex}

The notions in the next definition extend to graded rings. For the Noetherian
case, see~\cite{Matsumura}.
\begin{defn}\label{defn:RegularRings}
Let $A$ be a commutative ring.
\begin{itemize}
\item
$A$ is a \emph{regular local ring} if it is a local ring whose maximal ideal
$\mathfrak{m}$ is generated by a sequence $u_1,u_2,\ldots,u_n$, where~$n$ is
the Krull dimension of $A$; such a sequence is regular.
\item
$A$ is a \emph{regular ring} if for every maximal ideal $\mathfrak{m}\ideal A$,
the localization $A_{\mathfrak{m}}$ is a regular local ring.
\item
$A$ is a \emph{non-Noetherian regular local ring} if its maximal ideal is
generated by an infinite countable regular sequence $u_1,u_2,\ldots$.
\item
$A$ is a \emph{non-Noetherian regular ring} if for every maximal ideal
$\mathfrak{m}\ideal A$, the localization $A_{\mathfrak{m}}$ is a (possibly
non-Noetherian) regular local ring.
\end{itemize}
\end{defn}
\begin{prop}\label{prop:RegRings}
If $R_*$ is a \emph{(}possibly non-Noetherian\emph{)} regular local ring
with maximal ideal $\mathfrak{m}\ideal R_*$, then there is an $R$-module
$\kappa(\mathfrak{m})$ for which $\kappa(\mathfrak{m})_*=R_*/\mathfrak{m}$.
Therefore $\mathfrak{m}$ satisfies
\emph{Condition~(A)}.
\end{prop}
\begin{proof}
Given a (possibly infinite) regular sequence $u_1,u_2,\ldots$ which generates
$\mathfrak{m}$, we may follow the approach of~\cite[section~V.1]{EKMM} to
construct an $R$-module $\kappa(\mathfrak{m})$ which realizes $R_*/\mathfrak{m}$
as $\kappa(\mathfrak{m})_*=\pi_*\kappa(\mathfrak{m})$.
\end{proof}
\begin{cor}\label{cor:RegRings}
If $R_*$ is a regular ring then every maximal ideal $\mathfrak{m}\ideal R_*$
satisfies
\emph{Condition~(A)}.
\end{cor}

There is an associated \emph{Koszul complex}
\[
K_{*,*}=\Lambda_{R_*}(e_i:i\geq 1),
\]
with $e_i$ in bidegree $(1,|u_i|)$ and which is a differential graded
algebra with differential~$d$ given by
\[
d e_i=u_i.
\]
This provides a free resolution
\[
K_{\bullet,*}\lra R_*/\mathfrak{m}\ra0.
\]
Following the method of construction of the K\"unneth spectral sequence
in~\cite[section~IV.5]{EKMM}, we can use this to define a cell structure
on the $R$-module $\kappa(\mathfrak{m})$, with the cells corresponding
to the distinct monomials in the $e_i$.

Later, we will need multiplicative structures on our residue fields. These are 
provided by Angeltveit's result \cite[Theorem 4.2]{Angeltveit:E/I}.
\begin{thm} \label{thm:angeltveit}
If $R$ is a commutative $\mathbb{S}$-algebra whose coefficients are 
concentrated in even degrees and if an ideal $I \ideal R_*$ is generated by a 
regular sequence, then there is an $\mathbb{S}$-algebra structure on $R/I$ 
and $R \ra R/I$ is central. 
\end{thm}

\section{A finiteness result}\label{sec:Invertible}

General treatments of invertible objects in derived categories and
Picard groups may be found in~\cite{Fausk,FLM,JPM:PicardGps}. It is
standard that an invertible object in $\mathscr{D}_{R}$ is strongly
dualizable~\cite{HPS}. The following result on strongly dualizable
objects in $\mathscr{D}_{R}$ is taken from~\cite[proposition~2.1]{FLM}
(also see~\cite{JPM:EqtHtpyCohThy,JPM:PicardGps}).
\begin{prop}\label{prop:StrongDual-Cell2}
Let $X$ be an $R$-module. Then $X$ is strongly dualizable in $\mathscr{D}_{R}$
if and only if it is weakly equivalent to a retract of a finite cell $R$-module.
\end{prop}

Let $U$ be an \emph{invertible} $R$-module, \ie, $U$ is a cofibrant $R$-module 
and there is a cofibrant $R$-module~$V$ for which $U\wedge_R V\simeq R$. Then 
$V$ is a strong dual for $U$ and by
Proposition~\ref{prop:StrongDual-Cell2}, $U$ and $V$ are retracts of finite
cell $R$-modules.

The following Lemma allows us to apply coherence conditions to
topological settings.
\begin{lem}\label{lem:Coherent-cell}
Let $R$ be a commutative $\mathbb{S}$-algebra with coherent coefficient 
ring $R_*$.
\begin{enumerate}
\item[(a)]
Any finite cell $R$-module $M$ gives rise to a finitely generated
coherent $R_*$-module $M_*$.
\item[(b)]
Every retract $N$ of a finite $R$-cell module $M$ has finitely generated
coherent coefficients $N_*$.
\end{enumerate}
\end{lem}
\begin{proof}
As $M$ is built in finitely many steps via cofibre sequences of the form
\[
 \Sigma^{n} R \lra X \lra Y,
\]
its coefficients $M_*$ are built up out of exact couples of the form
\[
\xymatrix{ X_* \ar[rr] & & Y_* \ar[dl]\\
&  \Sigma^n R_* \ar[ul] &}
\]
Applying~\cite[theorem~3.1]{JCohen:Coherent}, we see that $M_*$ is
finitely generated coherent.

As retracts of cell $R$-modules correspond to finitely generated
submodules of $R_*$-modules $M_*$ as above, the claim follows.
\end{proof}
\begin{cor}\label{cor:Invtble->fg}
Suppose that $R_*$ is coherent and $U$ is an invertible $R$-module.
Then $U_*$ is a coherent $R_*$-module and hence it has a resolution
by finitely generated projective $R_*$-modules.
\end{cor}

\section{The connective case}\label{sec:Main}

Let $R$ be a connective commutative $\mathbb{S}$-algebra and let
$\mathfrak{m}\ideal R_*$ be a maximal ideal. Recall from~\cite{EKMM}
that there is an Eilenberg-Mac~Lane object $H(R_0/\mathfrak{m}_0)$
which is also a commutative $R$-algebra. We will view this as a
residue field $\kappa(\mathfrak{m})$. 
\begin{lem}\label{lem:MinRes-SS}
Let $R$ be a connective commutative $\mathbb{S}$-algebra and let
$\mathfrak{m}\ideal R_*$ be a maximal ideal for which
$(R_*)_{\mathfrak{m}}$ is coherent.
If $M,N$ are $R$-modules with $(N_*)_{\mathfrak{m}}$ coherent as
an $(R_*)_{\mathfrak{m}}$-module, then the $\mathrm{E}^2$-term of
the K\"unneth spectral sequence has the form
\begin{equation}\label{eqn:MinRes-SS}
\mathrm{E}^2_{p,*}=
\Tor^{{(R_\mathfrak{m})}*}_{p,*}(\kappa(\mathfrak{m})^R_*M_{\mathfrak{m}},
{(N_\mathfrak{m})}_*)
\cong \kappa(\mathfrak{m})^R_*M_{\mathfrak{m}}\oTimes_{R_*}Q_{p,*},
\end{equation}
where $Q_{\bullet, *}$ is a minimal resolution of ${(N_{\mathfrak{m}})}_*$.
\end{lem}
\begin{proof}
As $\pi_*(\kappa(\mathfrak{m})\wedge_R M)$ is an $(R_{\mathfrak{m}})_*$-module,
we can replace $R$ by $R_{\mathfrak{m}}$ and $N$ by its localization
$N_{\mathfrak{m}}$. Thus we might as well assume that $R_*$ is a coherent
local ring for the remainder of this proof.

We begin by choosing a free resolution $Q_{\bullet,*}\lra N_*\ra0$ of $N_*$
by finitely generated free $R_*$-modules. Using
Proposition~\ref{prop:Coherent-MinResn}, we can arrange this to be minimal.

Following~\cite{EKMM}, the $\mathrm{E}^2$-term for the K\"unneth spectral
sequence~\eqref{eqn:MinRes-SS} can be constructed using the above resolution,
giving
\[
\mathrm{E}^2_{p,*}=\Tor^{R_*}_{p,q}(\kappa(\mathfrak{m})^R_*M,N_*)=
\mathrm{H}_p(\kappa(\mathfrak{m})^R_*M\oTimes_{R_*}Q_{\bullet,*},\id\otimes
d_\bullet).
\]
By minimality this yields
\begin{equation}\label{eqn:MinRes-SS->Tor}
\mathrm{E}^2_{p,*}=\kappa(\mathfrak{m})^R_*M\oTimes_{R_*}Q_{p,*}.
\qedhere
\end{equation}
\end{proof}

We begin with a local result. Recall that a finitely generated module
over a local ring is projective if and only if it is free.
\begin{prop}\label{prop:Invtble-local}
Suppose that $R$ is connective and that $R_*$ is coherent and local
with maximal ideal $\mathfrak{m}\ideal R_*$. If $U$ is an invertible
$R$-module, then for some $k\in\Z$, $U\simeq\Sigma^k R$ and $U_*$ is
an invertible graded $R_*$-module.
\end{prop}
\begin{proof}
We follow the ideas and notation in the proof of Lemma~\ref{lem:MinRes-SS},
taking $M=U$ and $N=V$ where $U\wedge_R V\simeq R$.

In order to shorten notation we will write $\kappa=H(R_0/\mathfrak{m}_0)$.
By Corollary~\ref{cor:Invtble->fg} we can choose a minimal free resolution
$Q_{\bullet,*}\lra V_*\ra0$. Using~\eqref{eqn:MinRes-SS->Tor} we get as
$\mathrm{E}^2$-term of the K\"unneth spectral sequence
\[
\mathrm{E}^2_{p,*}=\Tor^{R_*}_{p,*}(\kappa^R_*U,V_*) \cong
\kappa^R_*U \otimes Q_{p,*}.
\]
Without loss of generality we can assume that $U$ is connective and $\pi_0(U) 
\neq 0$. The K\"unneth isomorphism
\[
\kappa_*^R(U) \otimes_{\kappa_*} \kappa_*^R(V) \cong \kappa_*^R(R)\cong\kappa_*
\]
forces $\kappa_*^R U$ to be free of rank one over
$\kappa_* \cong R_0/\mathfrak{m}_0$. Thus it must be concentrated
in degree zero and has to be isomorphic to
$R_0/\mathfrak{m}_0 = \kappa_0$. The whole spectral sequence is
concentrated in the first quadrant. No differential can hit the
entry in the $(0,0)$-coordinate. Therefore
\[
(\kappa_* \otimes_{R_*} Q_{0,*})_0 = \kappa_0 \otimes_{R_0} Q_{0,0}
\cong \kappa_0.
\]
As $Q_{0,0}$ is the zeroth homotopy group of some sum of $R$-spheres,
this forces $Q_{0,0}$ to be equal to $R_0$, in particular it is free
over $R_0$. The minimality of the resolution ensures that $Q_{p,0}$
is zero for $p>0$. Therefore the zero-line $\mathrm{E}^2_{p,0}$
vanishes except for $\kappa_0$ at $p=0$.

Inductively we assume that $Q_{p,i} = 0$ for all $p > 0, i \leq n$, and
that $Q_{0,i} \cong R_i$ for all $i \leq n$. Then the $(0,n+1)$-entry
in the $\mathrm{E}^2$-term cannot be hit by any differential, therefore
it must be an infinite cycle. As nothing else survives in total degree
$n+1$, we obtain
\begin{equation} \label{eqn:linearpart}
(\kappa_* \otimes_{R_*} Q_{0,*})_{n+1} = \kappa_{0} \otimes Q_{0,n+1}=0.
\end{equation}
This means that the whole module $Q_{0,n+1}$ gets killed by the relations
in the tensor product over $R_*$. We know that $R_{n+1}\subseteq Q_{0,n+1}$
because $Q_{\bullet,*}$ is an $R_*$-resolution and $Q_{0,0}\cong R_0$.
From~\eqref{eqn:linearpart} we know that there cannot be more in $Q_{0,n+1}$.

As $Q_{0,n+1} \cong R_{n+1}$, minimality ensures again that $Q_{p,n+1}= 0$
for all $p> 0$ and the induction is carried on.

Hence we obtain, that $Q_{0,q} \cong R_q$ for all $q \geq 0$ and the
higher terms in the resolution $Q_{p,q} \cong 0$ for $p >0$. This gives
$V_* \cong R_*$, proving the claim.
\end{proof}

We now use our local information to obtain a global result.
\begin{thm}\label{thm:Invtble-Main}
If $R$ is a connective commutative $\mathbb{S}$-algebra, such that
for every maximal ideal $\mathfrak{m} \ideal R_*$ the localization
is coherent. Then every invertible $R$-module spectrum $U$ has
invertible coefficients $U_*$.
\end{thm}
\begin{proof}
Let $V$ be an inverse for $U$. For each maximal ideal $\mathfrak{m}\ideal R_*$,
\[
\Tor^{R_*}_{s,*}(U_*,V_*)_{\mathfrak{m}}\iso
\Tor^{(R_*)_{\mathfrak{m}}}_{s,*}((U_*)_{\mathfrak{m}},(V_*)_{\mathfrak{m}})
\]
and also
\[
(U_{\mathfrak{m}}\wedge_{R_{\mathfrak{m}}} V_{\mathfrak{m}})
\simeq(U\wedge_R V)_{\mathfrak{m}} \simeq R_{\mathfrak{m}}.
\]
By our local result Proposition~\ref{prop:Invtble-local} we have
\[
U_{\mathfrak{m}}\simeq\Sigma^kR_{\mathfrak{m}}, \quad
V_{\mathfrak{m}}\simeq\Sigma^{-k}R_{\mathfrak{m}}.
\]
Hence for $s>0$,
\[
\Tor^{R_*}_{s,*}(U_*,V_*)_{\mathfrak{m}}=0.
\]

Now by a standard result on localizations~\cite{Matsumura},
\[
\Tor^{R_*}_{s,*}(U_*,V_*)=0\quad(s>0).
\]
So we find that the edge homomorphism $U_*\oTimes_{R_*}V_*\lra R_*$
of the K\"unneth spectral sequence is an isomorphism. Therefore $U_*$
is an invertible graded $R_*$-module with inverse $V_*$.
\end{proof}

\section{Eilenberg's condition}\label{sec:Connective}

For some important examples of spectra the coherence requirement
is too much to ask for. In~\cite{Eilenberg:syzygies}, Eilenberg
introduced conditions which ensure the existence of minimal
resolutions. We recall a particular case which then applies to
$\mathbb{S}$ and other commutative $\mathbb{S}$-algebras, leading
to important topological results.

Recall that a graded group $M_*$ is \emph{connective} if $M_n=0$
whenever $n<0$. Also, if $A_*$ is a connective graded commutative
local ring, then its unique maximal ideal $\mathfrak{m}\ideal A_*$
has components
\[
\mathfrak{m}_n=
\begin{cases}
\mathfrak{m}'& \text{if $n=0$}, \\
A_n& \text{otherwise},
\end{cases}
\]
where $A_0$ is local with maximal ideal $\mathfrak{m}'\ideal A_0$.
\begin{prop}
\label{prop:Eilenberg}
Let $A_*$ be a connective graded commutative local ring for
which $A_0$ is Noetherian and each $A_n$ is a finitely generated
$A_0$-module. Then every finitely generated $A_*$-module admits
a minimal resolution by free $A_*$-modules.
\end{prop}
\begin{proof}
See~\cite[proposition~14]{Eilenberg:syzygies}.
\end{proof}
\begin{ex}\label{ex:Eilenberg-S*}
The maximal ideals in the graded ring $\mathbb{S}_*$ have the
form
\[
\mathfrak{m}(p)_n=
\begin{cases}
(p)\ideal\Z& \text{if $n=0$}, \\
\ph{a}\mathbb{S}_n& \text{otherwise},
\end{cases}
\]
for rational primes $p>0$. On localizing we obtain the graded
local rings $(\mathbb{S}_*)_{(p)}$ which satisfy the requirements
of Proposition~\ref{prop:Eilenberg}.
\end{ex}
\begin{lem}\label{lem:Eilenberg-R}
Let $R$ be a commutative $\mathbb{S}$-algebra with connective
homotopy ring $R_*$ and let $\mathfrak{m}\ideal R_*$ be a maximal
ideal. If $(R_{\mathfrak{m}})_*=(R_*)_{\mathfrak{m}}$ satisfies
the requirements of\/ {\rm Proposition~\ref{prop:Eilenberg}},
then for any retract $W$ of a finite cell $R_{\mathfrak{m}}$-module,
there is a minimal resolution of $W_*$ by free
$(R_{\mathfrak{m}})_*$-modules.
\end{lem}
\begin{proof}
For a finite cell module $W$, this involves a straightforward
inductive verification that $W_*$ is finitely generated. But any
retract of such a $W$ has the same property.
\end{proof}
\begin{thm}\label{thm:Eilenberg-R}
Let $R$ be a commutative $\mathbb{S}$-algebra whose homotopy ring
$R_*$ is connective. Suppose that for every maximal ideal
$\mathfrak{m}\ideal R_*$, the localization
$(R_{\mathfrak{m}})_*=(R_*)_{\mathfrak{m}}$ has $(R_0)_{\mathfrak{m}}$
Noetherian and each $(R_n)_{\mathfrak{m}}$ is a finitely generated
$(R_0)_{\mathfrak{m}}$-module. If $U$ and $V$ are invertible
$R$-modules for which $U\wedge_R V\simeq R$, then we have
\[
U_*\otimes_{R_*}V_*\iso R_*
\]
and so $U_*$ is an invertible graded $R_*$-module. In particular,
there is a $k\in\Z$ for which $U_m=0=V_n$ whenever $m<k$ and $n<-k$
and then
\[
U_{k}\otimes_{R_0}V_{-k}\iso R_0,
\]
so $U_{k}$ is an invertible $R_0$-module.
\end{thm}

Using this, we obtain the following well-known result of~\cite{NS:Picard,HMS}.
\begin{ex}\label{ex:Eilenberg-S}
Taking $R=\mathbb{S}$ and recalling Example~\ref{ex:Eilenberg-S*},
we see that $U\wedge V\simeq\mathbb{S}$ implies that for some $k\in\Z$
as in the Theorem, $U_k$ is an invertible $\Z$-module and
$U_k\iso\Z\iso V_{-k}$, hence $U_{*+k}\iso\mathbb{S}_*\iso V_{*-k}$.
It follows that $\Sigma^{-k}U \simeq\mathbb{S}\simeq\Sigma^k V$.
\end{ex}

Other examples include $MU$, $MSp$, $ku$, $ko$ and the connective
spectrum of topological modular forms $\mathit{tmf}$. The last
example is known to be a commutative $\mathbb{S}$-algebra and its
homotopy ring is computed in~\cite{Bauer:tmf,Rezk:tmf}; it has
$\pi_0\mathit{tmf}=\Z$ and satisfies the conditions of
Proposition~\ref{prop:Eilenberg}.

So far we did not give any proof in the case of Eilenberg-Mac~Lane
spectra over arbitrary commutative rings. For the sake of completeness
we add this result here, although it is probably well-known.
\begin{prop}\label{prop:EMSpectra}
Let $A$ be a commutative ring with unit. Then for every invertible
$HA$-module spectrum $U$, $U_*$ is an invertible graded $A$-module.
\end{prop}
\begin{proof}
The proof we give here is an elementary adaption of Fausk's proof
in~\cite[3.2,3.3]{Fausk} to our setting. Without loss of generality
we can assume that $U$ has its first non-vanishing homotopy group in
degree zero.

Assume first that $A$ is a local ring. Let $V$ be an inverse of $U$
over $HA$. We know that $U_0 \otimes_A V_0 \cong A$ because nothing
else can hit the zeroth homotopy group in
\[
\mathrm{E}^2_{p,q} = \Tor_{p,q}^A(U_*,V_*) \Lra \pi_0(HA) = A.
\]
As $A$ is local, the only invertible $A$-modules are the ones which
are isomorphic to $A$. In particular $U_0$ and $V_0$ are free.
Therefore the rest of the $(p,0)$-line vanishes. This forces the
$(1,0)$-entry to survive, so it has to be trivial, which means that
$U_1=0=V_1$. Iteratively, we can clear out the whole $\mathrm{E}^2$-page
except for the $(0,0)$-entry. In particular, for all $p>0$,
\[
\Tor_{p,q}^A(U_*,V_*) = 0.
\]
A local-to-global argument then proves the result in general.
\end{proof}

\section{Small global dimension}\label{sec:Invtblenon-Conn}


The crucial point in our proofs is the collapsing of
the K\"unneth spectral sequence. For commutative $\mathbb{S}$-algebras
with residue fields we gain an analog of Theorem~\ref{thm:Invtble-Main}
as long as we can exclude non-trivial differentials.
\begin{thm}\label{thm:nonconnMain}
Let $R$ be a commutative $\mathbb{S}$-algebra, such that for every
maximal ideal $\mathfrak{m} \ideal R_*$ the ring ${(R_*)}_\mathfrak{m}$
is coherent and assume that $R$ satisfies\/ {\rm Condition~(A)}
and has a structure of a ring spectrum on each of its residue fields. If $R_*$
has global dimension at most~$2$ then every invertible $R$-module
spectrum $U$ has invertible coefficients $U_*$.
\end{thm}
\begin{proof}
As we can perform a local-to-global argument, we may as well assume
that $R_*$ is local and coherent. Let $V$ be an inverse of $U$. The
existence of a residue field $\kappa$ which is a ring spectrum ensures
that $\kappa_*^R(U)$ is a $\kappa_*$-vector space. The K\"unneth map
\[
\kappa^R_*(U) \otimes_{\kappa_*} \kappa^R_*(V) \lra \kappa_*
\]
has to be an isomorphism. Therefore we can set $\kappa^R_*(U)\cong\kappa_*$.

Together with the existence of a minimal resolution of $V_*$ this guarantees
that the $\mathrm{E}^2$-term of the K\"unneth spectral sequence is given by
\[
\mathrm{E}^2_{p,*}=\Tor_{p,*}^{R_*}(\kappa_*, V_*)
\cong \kappa_*\otimes_{R_*} Q_{p,*}.
\]
If the global dimension of $R_*$ is at most $1$, this spectral sequence
is concentrated in two columns. Therefore there cannot be any non-trivial
differentials. The abutment of the spectral sequence is $\kappa_*$. As
$Q_{p,*}$ is a resolution of $V_*$, $Q_{0,*}$ cannot be trivial. If
$Q_{1,*}$ were non-trivial a dimension count leads to a contradiction.
Similarly, $Q_{0,*}$ must be free of rank one over $R_*$ and therefore
$V_* \cong \Sigma^kR_*$ for some $k \in \Z$.

In the case of global dimension $2$ the $\mathrm{E}^2$-term of the
K\"unneth spectral sequence converging to $\kappa_*$
\[
\mathrm{Tor}_{p,*}^{R_*}(\kappa_*, V_*) \cong \kappa_*
\otimes_{R_*} Q_{p,*}
\]
has only three non-trivial columns
\[
\xymatrix{
{} &{} & {}  \\
{} &{} &{}  \\
{\kappa_* \otimes_{R_*} Q_{0,*}} \ar@{.}[d] \ar@{.}[uu]&
{\kappa_* \otimes_{R_*} Q_{1,*}} \ar@{.}[d] \ar@{.}[uu]&
{\kappa_* \otimes_{R_*} Q_{2,*}} \ar[llu]_{\phantom{bla}d^2} \ar@{.}[d] 
\ar@{.}[uu] \\
{} &{} & {}
}
\]
There are two possible cases. If the differential $d^2$ is trivial,
then the $\mathrm{E}^2$-term is the $\mathrm{E}^\infty$-term. As
the spectral sequence has a one-dimensional abutment, we see as
before that $V_* \iso \Sigma^k R_*$ for some $k\in\Z$. On the other
hand, if $d^2$ is non-trivial, then we have a non-trivial map between
the $0$-column and the $2$-column. But the entries in the $1$-column
are infinite cycles. So either they are trivial or at most one-dimensional.
If they are trivial then the $2$-column has to be trivial as well and
we get a contradiction. If they are non-trivial, then we can conclude
that the $d^2$-differential must be an isomorphism and that
$Q_{1,*} \iso \Sigma^\ell R_*$ for some $\ell\in \Z$. Therefore up
to suspensions the resolution $Q_{\bullet,*}$ is of the form
\[
R_*^n \lla R_* \lla R_*^n.
\]
If $n$ were bigger than one, then this could not be a resolution. For
$n=1$ the differentials must be given by multiplication by some element
in $R_*$. Such maps have a non-trivial kernel, therefore this gives no
resolution either.
\end{proof}

We can loosen the requirements on $R$ a little bit by referring to
Condition~(B).
\begin{prop}\label{prop:CondB-case}
Let $R$ be a commutative $\mathbb{S}$-algebra such that for every
maximal ideal $\mathfrak{m} \ideal R_*$ the ring ${(R_*)}_\mathfrak{m}$
is coherent and satisfies\/ {\rm Condition~(B)} with residue fields
which are ring spectra. If $R_*$ has global dimension at most~$2$ then
for every invertible $R$-module spectrum $U$, $U_*$ is an invertible
graded $R_*$-module.
\end{prop}
\begin{proof}
Let $R\lra R'$ be the unit of a suitable $R$-algebra as required in
Condition~(B) and let $V' = R' \wedge_R V$. Coherence of $R$ guarantees
the existence of a minimal resolution $Q_{\bullet,*}\lra V_*\ra0$ of
$V_*$. Flatness of $R'$ over $R$ ensures that
\[
Q'_{\bullet,*}\lra V'_* = R'_* \otimes_{R_*} V_* \ra0
\]
is still a resolution and as we assumed the map $R_* \lra R'_*$ to
be local, this resolution is minimal. Using the proof of
Theorem~\ref{thm:nonconnMain} and an argumentation as in the proof
of Proposition~\ref{prop:Invtble-local} for $R'$ and $U' = R'\wedge_RU$,
we see that $Q'_{p,*}=0$ when $p>0$. As $Q'_{p,*}=R'_*\oTimes_{R_*}Q_{p,*}$
and $Q_{p,*}$ is free over $R$, we must have $Q_{p,*}=0$ for $p>0$.
Thus $Q_{0,*}\cong V_*$ and the K\"unneth spectral sequence
\[
\mathrm{E}^2_{p,q}=\Tor^{R_*}_{p,q}(U_*,V_*)
                             \Lra\pi_{p+q}(U\wedge_R V)=R_{p+q}
\]
collapses and the edge homomorphism $U_*\oTimes_{R_*}V_*\lra R_*$
is an isomorphism, so $U_*$ is invertible with inverse $V_*$.
\end{proof}

%
\begin{ex}\label{ex:dim<=2}
Theorem~\ref{thm:nonconnMain} and Proposition~\ref{prop:CondB-case}
cover the examples of the first two Lubin-Tate spectra $E_1$ and
$E_2$ and their close relatives the completed Johnson-Wilson spectra
$\widehat{E(1)}$ and $\widehat{E(2)}$, as well as the Adams summand
$E(1)$. Complex periodic $K$-theory and real periodic $K$-theory with
$2$-inverted, $KO[1/2]$, fulfills the requirements as well.
\end{ex}

\section{Noetherian complete local regular rings}\label{sec:RegRings}

In the following we extend the results of Section~\ref{sec:Invtblenon-Conn}
to commutative $\mathbb{S}$-algebras whose coefficients have higher
global dimension. However, we have to impose regularity conditions.
The method of proof is adapted from that 
of~\cite[theorem~1.3, pp.117,118]{HMS}.

We will make use of the algebraic theory of Noetherian regular rings
and their finite modules for which we refer to~\cite{Bruns&Herzog,Matsumura}.
We begin with some local results.
\begin{assump}\label{assump:RegRings}
Throughout this section, $R$ will be a commutative $\mathbb{S}$-algebra
for which $R_* = R_0[u, u^{-1}]$ with $|u|=2$. We assume that $R_0$ is
a complete Noetherian local regular ring whose maximal ideal
$\mathfrak{m}\ideal R_0$ is generated by a regular sequence $u_1,\ldots,u_n$,
where $n$ is the Krull dimension of $R_*$. We could view $R_*$
and its modules as $\Z/2$-graded $R_0$-modules.
\end{assump}

Theorem~\ref{thm:angeltveit} then applies. 
\begin{lem}\label{lem:Angeltveit-E/I}
For each prime ideal $\mathfrak{p}\ideal R_*$, there is an $R$-algebra
realizing the $R_*$-algebra $R_*/\mathfrak{p}$. Hence the graded
residue field
\[
\kappa(\mathfrak{p})_*=(R_*/\mathfrak{p})_{\mathfrak{p}}
                 =(R_*)_{\mathfrak{p}}/(R_*)_{\mathfrak{p}}\mathfrak{p}
\]
can be realized as an $R_{\mathfrak{p}}$-algebra and so $R_{\mathfrak{p}}$
has a residue field $\kappa(\mathfrak{p})$.
\end{lem}

For $R$-modules $M$ and $N$, $\kappa(\mathfrak{m}) \wedge_R M$ and
$\kappa(\mathfrak{m}) \wedge_R N$ are $\kappa(\mathfrak{m})$-left
modules, and we can consider $\kappa(\mathfrak{m}) \wedge_R M$ as
a right $\kappa(\mathfrak{m})$-module spectrum via the action of
$\kappa(\mathfrak{m})$ on itself by right multiplication. Since
$R$ is central in $\kappa(\mathfrak{m})$, this is well-defined. 
\begin{cor}\label{cor:Angeltveit-E/I}
Let $M$ and $N$ be $R$-modules. Then there is a K\"unneth isomorphism
\[
\kappa(\mathfrak{m})^R_*(M)\otimes_{\kappa(\mathfrak{m})_*}
                      \kappa(\mathfrak{m})^R_*(N)
\iso\kappa(\mathfrak{m})^R_*(M\wedge_R N).
\]
If\/ $U$ is an invertible $R$-module then
$\dim_{\kappa(\mathfrak{m})_*}\kappa(\mathfrak{m})^R_*(U)=1$.
\end{cor}
\begin{proof}
In our case, the K\"unneth spectral sequence of~\cite[theorem IV.4.1]{EKMM}
\[
\mathrm{E}^2_{p,q}=
\Tor^{\kappa(\mathfrak{m})_*}_{p,q}(\kappa(\mathfrak{m})^R_*(M),
\kappa(\mathfrak{m})^R_*(N))
\Lra
\kappa(\mathfrak{m})^R_{p,q}(M\wedge_R N),
\]
collapses because $\kappa(\mathfrak{m})_*$ is a graded field.
When $U\wedge_R V\simeq R$, we have
\[
\kappa(\mathfrak{m})^R_*(U)\otimes_{\kappa(\mathfrak{m})_*}
\kappa(\mathfrak{m})^R_*(V)\iso\kappa(\mathfrak{m})_*,
\]
hence
\[
\dim_{\kappa(\mathfrak{m})_*} \kappa(\mathfrak{m})^R_*(U)=1
=\dim_{\kappa(\mathfrak{m})_*} \kappa(\mathfrak{m})^R_*(V).
\qedhere
\]
\end{proof}

We will need some technical results about killing regular sequences
in~$R$. We make use of the results of~\cite[lemma~V.1.5]{EKMM}.
\begin{lem}\label{lem:cofseq}
For every sequence $(u_1^{i_1},\ldots, u_n^{i_n})$ with $i_j > 1$,
there are cofibre sequences
\begin{equation}\label{eqn:cof1}
R/(u_1^{i_1},\ldots, u_n^{i_n}) \xrightarrow{u_j}
R/(u_1^{i_1},\ldots, u_n^{i_n}) \xrightarrow{\ph{u_j}}
R/(u_1^{i_1},\ldots,u_j^1,\ldots, u_n^{i_n})
\vee
\Sigma R/(u_1^{i_1},\ldots,u_j^1,\ldots, u_n^{i_n})
\end{equation}
and
\begin{equation} \label{eqn:cof2}
R/(u_1^{i_1},\ldots,u_j^{i_j-1}, \ldots, u_n^{i_n})
\lra
R/(u_1^{i_1},\ldots,u_j^{i_j}, \ldots, u_n^{i_n})
\lra
R/(u_1^{i_1},\ldots,u_j^1, \ldots, u_n^{i_n}).
\end{equation}
\end{lem}
\begin{proof}
The cofibre of the multiplication map by $u_j$ on
$R/(u_1^{i_1},\ldots,u_n^{i_n})$ can be identified as follows. As
the variables behave independently we might just consider the case
of one $u_j$. Then we get the following diagram of cofibre sequences.
\[
\xymatrix{
{R} \ar[r]^{u^i_j} \ar[d]_{u_j} & {R} \ar[d]^{u_j} \ar[r] & {R/u_j^i}
\ar[d]^{u_j} \\
{R} \ar[r]^{u_j^i} \ar[d] & {R} \ar[d] \ar[r] & {R/u^i_j} \ar[d] \\
{R/u_j} \ar[r]^{u_j^i} & {R/u_j} \ar[r] & {\cofibre(u_j)}
}
\]
As multiplication by $u^i_j$ is nullhomotopic on $R/u_j$, the cofibre
of the multiplication by $u_j$ splits as $R/u_j \vee \Sigma R/u_j$.

\noindent
For the second sequence, consider
\begin{multline*}
R/(u_1^{i_1},\ldots, u_{j-1}^{i_{j-1}},u_{j+1}^{i_{j+1}}, \ldots,u_n^{i_n})
\xrightarrow{u_j^{i_j-1}}
R/(u_1^{i_1},\ldots, u_{j-1}^{i_{j-1}}, u_{j+1}^{i_{j+1}}, \ldots,u_n^{i_n}) \\
\xrightarrow{\ph{u_j^{k-1}}}
R/(u_1^{i_1},\ldots, u_{j-1}^{i_{j-1}},u_j^{i_j-1},u_{j+1}^{i_{j+1}},\ldots,u_n^{i_n}).
\end{multline*}
There is a canonical projection map
\[
R/(u_1^{i_1},\ldots, u_{j-1}^{i_{j-1}}, u_{j+1}^{i_{j+1}}, \ldots,u_n^{i_n})
\lra R/(u_1^{i_1},\ldots,u_j^{i_j}, \ldots, u_n^{i_n})
\]
which we can compose with multiplication by $u_j$. When precomposed
with multiplication by $u_j^{i_j-1}$, this map becomes nullhomotopic,
thus by~\cite[lemma~V.1.5]{EKMM} it factors through a map
\[
R/(u_1^{i_1},\ldots,u_j^{i_j-1}, \ldots, u_n^{i_n})
\lra
R/(u_1^{i_1},\ldots,u_j^{i_j}, \ldots, u_n^{i_n})
\]
whose cofibre is easily identified with
$R/(u_1^{i_1},\ldots,u_j^1,\ldots, u_n^{i_n})$.
\end{proof}
\begin{lem}\label{lem:K(m)*U}
If\/ $U$ is an invertible $R$-module spectrum, then for all sequences
$(u_1^{i_1}, \ldots, u_n^{i_n})$ with $\sum_{k=1}^n i_k \geq n$ and each 
$i_k \geq 1$, $(R/(u_1^{i_1}, \ldots, u_n^{i_n}))^R_*(U)$ is a cyclic
$R_*$-module.
\end{lem}
\begin{proof}
By suspending $U$ if necessary and appealing to Lemma~\ref{cor:Angeltveit-E/I},
we may as well assume that $\kappa(\mathfrak{m})^R_*(U)\cong \kappa(\mathfrak{m})_*$
is concentrated in even degrees.
We prove the claim by induction on $m=\sum_{k=1}^n i_k$. For $m = n$
the result is clear since
\[
(R/(u_1^{i_1}, \ldots, u_n^{i_n}))^R_*(U)=\kappa(\mathfrak{m})^R_*(U) \cong \kappa(\mathfrak{m})_*.
\]
Now let $m>n$ and assume that the result for all sequences of the above
form with $ m>\sum_{k=1}^n i_k$. Using the cofibre sequence~\eqref{eqn:cof2},
we see that the module $R/(u_1^{i_1}, \ldots,u_n^{i_n}) \wedge_R U$ has
homotopy groups which are concentrated in even degrees. From the cofibre
sequence~\eqref{eqn:cof1} we can read off that multiplication by $u_j$
on $(R/(u_1^{i_1}, \ldots, u_n^{i_n}))^R_*(U)$ has quotient
$(R/(u_1^{i_1},\ldots,u_j^1,\ldots, u_n^{i_n}))^R_*(U)$ which is cyclic
by assumption.

Notice that all three terms are finitely generated $R_*$-modules and
the image of the multiplication by $u_j$ is contained in the submodule
generated by the maximal ideal of $R_*$, hence by Nakayama's Lemma,
the module $(R/(u_1^{i_1},\ldots,u_j^{i_j}, \ldots, u_n^{i_n}))^R_*(U)$
is cyclic.
\end{proof}
\begin{lem}\label{lem:consistent}
If\/ $U$ is an invertible $R$-module spectrum, then for every sequence
$(u_1^{i_1}, \ldots, u_n^{i_n})$ with $\sum_{k=1}^n i_k \geq n$ and
$i_k \geq 1$, up to suspension, there is an isomorphism 
\[
(R/(u_1^{i_1},\ldots,u_n^{i_n}))^R_*(U)
                                \cong R_*/(u_1^{i_1},\ldots,u_n^{i_n}).
\]
Furthermore, these isomorphisms are compatible with the projection maps
\[
(R/(u_1^{i_1}, \ldots,u_j^{i_j + 1},\ldots, u_n^{i_n}))^R_*(U)
\lra
(R/(u_1^{i_1}, \ldots,u_j^{i_j},\ldots, u_n^{i_n}))^R_*(U).
\]
\end{lem}
\begin{proof}
There is a canonical cofibre sequence
\[
R/(u_1^{i_1},\ldots, u_j, \ldots, u_n^{i_n}) \lra
R/(u_1^{i_1}, \ldots,u_j^{i_j + 1},\ldots, u_n^{i_n}) \lra
R/(u_1^{i_1}, \ldots,u_j^{i_j},\ldots, u_n^{i_n}).
\]
As everything in sight is concentrated in even degrees, for each
$\ell\in\Z$ we get the two short exact sequences
\begin{multline*}
0\ra
R_{2\ell}/(u_1^{i_1},\ldots, u_j, \ldots, u_n^{i_n}) \lra
R_{2\ell}/(u_1^{i_1}, \ldots,u_j^{i_j + 1}, \ldots, u_n^{i_n}) \\
\lra R_{2\ell}/(u_1^{i_1}, \ldots, u_j^{i_j},\ldots,u_n^{i_n})
\ra0
\end{multline*}
and
\begin{multline*}
0\ra
(R/(u_1^{i_1},\ldots, u_j, \ldots, u_n^{i_n}))^R_{2\ell}(U) \lra
(R/(u_1^{i_1}, \ldots, u_j^{i_j + 1},\ldots, u_n^{i_n}))^R_{2\ell}(U) \\
\lra (R/(u_1^{i_1}, \ldots, u_j^{i_j},\ldots, u_n^{i_n}))^R_{2\ell}(U)
\ra0.
\end{multline*}
We start with the isomorphism $\kappa(\mathfrak{m})_*^R(U)\iso\kappa(\mathfrak{m})_*$.
As every $R/(u_1^{i_1}, \ldots,u_j^{i_j},\ldots, u_n^{i_n}))^R_*(U)$
is cyclic we can choose
epimorphisms
\[
R_{2\ell}/(u_1^{i_1}, \ldots,u_j^{i_j},\ldots, u_n^{i_n})
\lra
(R/(u_1^{i_1}, \ldots,u_j^{i_j},\ldots,u_n^{i_n}))^R_{2\ell}(U)
\]
which make the following diagram commute.
\[
\xymatrix@C=0.6cm@R=0.7cm{
{0} \ar[r]
&
{\scriptstyle R_{2\ell}/(u_1^{i_1},\ldots, u_j, \ldots,u_n^{i_n})} \ar[r] \ar[d] &
{\scriptstyle R_{2\ell}/(u_1^{i_1}, \ldots,u_j^{i_j +1}, \ldots, u_n^{i_n})} \ar[r] \ar[d] &
{\scriptstyle R_{2\ell}/(u_1^{i_1}, \ldots, u_j^{i_j},\ldots, u_n^{i_n})} \ar[r] \ar[d] &
{0}
\\
{0} \ar[r] &
{\scriptstyle (R/(u_1^{i_1},\ldots, u_j, \ldots,u_n^{i_n}))^R_{2\ell}(U)} \ar[r] &
{\scriptstyle (R/(u_1^{i_1}, \ldots, u_j^{i_j + 1},\ldots,u_n^{i_n}))^R_{2\ell}(U)}  \ar[r] &
{\scriptstyle (R/(u_1^{i_1}, \ldots, u_j^{i_j},\ldots,u_n^{i_n}))^R_{2\ell}(U)} \ar[r] & {0}}
\]
Now an induction over $m = \sum_{j=1}^n i_j$ proves the claim.
\end{proof}
\begin{thm}\label{thm:RegNoethLocal}
If $R$ satisfies  {\rm Assumption \ref{assump:RegRings}},
then every invertible $R$-module is equivalent to a suspension of~$R$.
\end{thm}
\begin{proof}
Again, we may suspend $U$ if necessary to ensure that
$\kappa(\mathfrak{m})_*^R(U)\cong\kappa(\mathfrak{m})_*$.


Lemma~\ref{lem:consistent} ensures that the identifications
\[
R_*/(u_1^{i_1}, \ldots,u_j^{i_j},\ldots, u_n^{i_n}) \cong
(R/(u_1^{i_1}, \ldots,u_j^{i_j},\ldots, u_n^{i_n}))_*^R(U)
\]
are consistent with the projection maps in the inverse system
for $\holim R/(u_1^{i_1}, u_2^{i_2}, \ldots, u_n^{i_n}) \wedge_R U$.
Since $R_*$ is a Noetherian complete local ring 
$\varprojlim R_*/(u_1^{i_1}, u_2^{i_2}, \ldots, u_n^{i_n}) = 
\varprojlim_{\ell} R_*/\mathfrak{m}^\ell$. As $U$ is a finite
cell $R$-module, we have
\[
\holim R/(u_1^{i_1}, u_2^{i_2}, \ldots, u_n^{i_n}) \wedge_R U
\simeq
R\wedge_R U \simeq U.
\]
Using the above description of
$(R/(u_1^{i_1}, u_2^{i_2}, \ldots, u_n^{i_n})^R_*(U)$ we find that
\[
\holim R/(u_1^{i_1}, u_2^{i_2}, \ldots, u_n^{i_n})\wedge_R U
\simeq
\holim R/(u_1^{i_1}, u_2^{i_2}, \ldots, u_n^{i_n}) \simeq R.
\]
Therefore we have $U\simeq R$. In the general case, $U$ might be
equivalent to $\Sigma R$.
\end{proof}

The proof of the following more general result involves a standard
local-to-global argument.
\begin{thm}\label{thm:CompleteRegular}
Let $R$ be a commutative $\mathbb{S}$-algebra such that the
localization of~$R$ at any maximal ideal $\mathfrak{m}\ideal R_*$
satisfies  {\rm Assumption \ref{assump:RegRings}}.
Then for every invertible $R$-module $U$, $U_*$ is an invertible
$R_*$-module.
\end{thm}
\begin{rem}\label{rem:afterCompleteRegular}
Assumption~\ref{assump:RegRings} is not optimal. For example, one
might loosen the requirement that $R_*$ is $2$-periodic and replace
this by a periodicity of degree $2\ell$ for some $\ell$. One might
also wish to allow that the generators of the maximal ideal then
lie in degrees different from zero. This is no problem if one takes
appropriate suspensions into account in the numerous cofibre sequences.
Last but not least there might be cases of infinite Krull dimension
that are tractable.
\end{rem}

As an application we discuss invertible modules over a group ring
$R[C_{p^\ell}]$. 

\begin{prop}\label{prop:GpRing}
Let $R$ be a commutative $\mathbb{S}$-algebra  which satisfies {\rm Assumption 
\ref{assump:RegRings}}. Suppose that $p$ is a prime which is not contained
in the maximal ideal $\mathfrak{m}$ and $R_0$ contains a primitive
$p^\ell$-th root of unity. Then for every invertible
$R[C_{p^\ell}]$-module $U$, $U_*$ is an invertible
$R_*[C_{p^\ell}]$-module. 
\end{prop}
\begin{proof}
Since $p$ is invertible and $R_0$ contains enough roots of unity, 
there is a complete set of orthogonal idempotents $e_i$ with
$i=1,\ldots,p^\ell$. These can be realized as elements of 
$\pi_0R[C_{p^\ell}]$ and therefore as self maps of $U$ whose images
$e_iU$ give a wedge decomposition  
$$U \sim \bigvee_{i=1}^{p^\ell} e_iU.$$
As argued in \cite[Theorem 2.3.2]{AB&BR:Galois} all summands $e_iU$ are 
invertible $R$-modules
and are therefore  equivalent to $R$ up to suspension. Thus 
$$ U \sim \bigvee_{i=1}^{p^\ell} R \sim R[C_{p^\ell}] $$
up to suspension. 
\end{proof}
This result extends to arbitrary finite abelian groups $G$ for which the 
primes  
dividing the order of the group are not contained in the maximal ideal 
$\mathfrak{m} \ideal R_0$, because there exists a more general decomposition 
into 
eigenspaces of characters $\chi \in \Hom(G,R_0^\times)$. 
\begin{ex}
Consider the Lubin-Tate spectrum $E_n$ with ${(E_n)}_0 = W\F_{p^n}[[u_1,
\ldots, u_{n-1}]]$. This ring is local, complete and regular. Let 
$W\F_{p^n}^{\mathrm{nr}}$ denote the  maximal unramified extension of 
$W\F_{p^n}$. There is a commutative $\mathbb{S}$-algebra $E_n^\mathrm{nr}$ 
with  coefficient ring 
$${E_n^\mathrm{nr}}_* = W\F_{p^n}^{\mathrm{nr}}[[u_1,
\ldots, u_{n-1}]][u^{\pm 1}] \quad \mathrm{with } \quad |u| = 2. $$
Every invertible $E_n$-module has invertible coefficients and as 
${E_n^\mathrm{nr}}_*$  
contains enough roots of unity every invertible module over 
$E_n^\mathrm{nr}[G]$ has 
invertible coefficients whenever $G$ is a finite abelian group whose order is 
not divisible by $p$. As ${E_n^\mathrm{nr}}_*[G]$ is semilocal, the Picard 
group 
$\Pic(E_n^\mathrm{nr}[G])$ is trivial. Using \cite{R:Opusmagnus}
or~\cite[Example~1.4.8]{AB&BR:Galois} this shows that every finite 
Galois extension of $E_n^\mathrm{nr}$ with Galois group as above possesses a 
normal basis. 
\end{ex}
For a  connective commutative $\mathbb{S}$-algebra $R$ with coherent 
coefficients and an arbitrary finite abelian groups $G$, the group ring 
$R_*[G]$ is coherent \cite[corollary 1.2]{Harris} and therefore in this case  
invertible $R[G]$-modules have invertible coefficients. 
\section{Examples}\label{sec:Examples}

We will now restate our earlier results in terms of these Picard groups.
\begin{thm}\label{thm:Pic}
For a commutative $\mathbb{S}$-algebra $R$, there is a monomorphism
of abelian groups
\[
\Phi\: \Pic(R_*) \lra \Pic(R).
\]
Furthermore, if $R$ satisfies the conditions of\/
{\rm Theorem~\ref{thm:Invtble-Main}}, {\rm Theorem~\ref{thm:Eilenberg-R}}
or of\/ {\rm Theorem~\ref{thm:CompleteRegular}}, then $\Phi$ is an
isomorphism.
\end{thm}

Thus $\Pic(R_*)\iso\Pic(R)$ in all of the following cases.
\begin{itemize}
\item
$R=H A$, where $A$ is a commutative ring.
\item
$R=MU/I$, where $I\ideal MU_*$ is a finitely generated
ideal for which
$MU/I$ is a commutative $\mathbb{S}$-algebra.
\item
$KU$, $KO[1/2]$, $ku$, and $ko$.
\item
$\mathit{tmf}$ at a prime $p$, $E(1)$,
$BP\<1\>$, $\widehat{E(1)}$, and $\widehat{E(2)}$.
See~\cite{AB&BR:GammaCoh-Numerical} for the existence of commutative
$\mathbb{S}$-algebra structures on some of these.
\item
$MSp$, $MSpin$ and $MSU$.
\item
$E_n$ for any $n$ and $p$. 
\end{itemize}

We close with a counterexample which originates from Galois theory of
commutative $\mathbb{S}$-algebras. In~\cite[Theorem~2.5.1]{AB&BR:Galois},
we show that every finite abelian Galois extension $B/A$ with Galois
group $G$ gives rise to an element in the Picard group of the group
ring $A[G]$.
\begin{ex}\label{ex:KU}
Complex periodic $K$-theory, $KU$, is a Galois extension of the real
periodic $K$-theory spectrum $KO$, whose Galois group is $C_2$, the
cyclic group of order~$2$ (see~\cite{R:Opusmagnus}
or~\cite[Example~1.4.8]{AB&BR:Galois}). Therefore we obtain
\[
KU \in \Pic(KO[C_2]).
\]
But $KU_*$ is not projective over $KO_*$, therefore it cannot
be projective over $KO_*[C_2]$. In particular, the coefficient
module $KU_*$ is not an element in $\Pic(KO_*[C_2])$.
\end{ex}

\end{document}